\documentclass{amsart}[10pt]

\usepackage{amscd,amssymb,latexsym,url,verbatim,graphicx,color}

\usepackage{cases,amsmath}

\usepackage{txfonts}

\usepackage{tikz,tikz-cd}
\usetikzlibrary{decorations.pathmorphing}

\usepackage[dvips]{epsfig}

\title
{Topology of scrambled simplices}

\author{Dmitry N. Kozlov}

\address{Department of Mathematics, University of Bremen, 28334
  Bremen, Federal Republic of Germany}

\email{dfk@math.uni-bremen.de}

\newtheorem{theorem}{Theorem}[section]
\newtheorem{df}[theorem]{Definition}
\newtheorem{thm}[theorem]{Theorem} 
\newtheorem{prop}[theorem]{Proposition}
\newtheorem{lm}[theorem]{Lemma}

\newtheorem{rem}[theorem]{Remark}

 \newcommand{\nin}{\noindent}
\newcommand{\pr}{\nin{\bf Proof.} }

\newcommand{\bd}{{\text{\rm bd}\,}}

\newcommand{\ce}{{\mathcal E}}
\newcommand{\cf}{{\mathcal F}}

\newcommand{\cp}{{\mathcal P}}

\newcommand{\zz}{{\mathbb Z}}

\newcommand{\da}{\Delta}
\newcommand{\lra}{\leftrightarrow}
\newcommand{\hra}{\hookrightarrow}

\newcommand{\ra}{\rightarrow}
\newcommand{\sm}{\setminus}
\newcommand{\supp}[1]{\text{\rm supp}\,(#1)}

\newcommand{\ti}{\tilde}
\newcommand{\wti}{\widetilde}

\newcommand{\id}{\textrm{id}}
\newcommand{\im}{\text{\rm Im}}

\newcommand{\dar}{\downarrow}
\newcommand{\succl}{\succ_{\text{\rm lex}}}
\newcommand{\succc}{\succ_{\text{\rm colex}}}
\newcommand{\alt}{\text{\rm alt}}
\newcommand{\sgn}{\text{\rm sgn}}

\numberwithin{equation}{section}
\numberwithin{figure}{section}
\numberwithin{table}{section}

\begin{document}

\begin{abstract}
In this paper we define a~family of topological spaces, which contains
and vastly generalizes the higher-dimensional Dunce hats.  Our
definition is purely combinatorial, and is phrased in terms of
identifications of boundary simplices of a~standard $d$-simplex. By
virtue of the construction, the obtained spaces may be indexed by
words, and they automatically carry the structure of a~$\da$-complex.

As our main result, we completely determine the homotopy type of these
spaces.  In fact, somewhat surprisingly, we are able to prove that
each of them is either contractible or homotopy equivalent to
an~odd-dimensional sphere. We develop the language to determine the
homotopy type directly from the combinatorics of the indexing word.

As added benefit of our investigation, we are able to emulate the
Dunce hat phenomenon, and to obtain a~large family of both
$\da$-complexes, as well as simplicial complexes, which are
contractible, but not collapsible.
\end{abstract}

\maketitle

\section{A combinatorial family of $\da$-complexes}

\subsection{Introduction} $\,$

\nin Imagine we are given a~$d$-simplex $\da^d$, viewed as
a~$\da$-complex\footnote{We use the terminology of Hatcher,
  see~\cite{Ha02}; alternatively these spaces were called {\it
    semisimplicial sets} in~\cite{EZ}, {\it triangulated spaces}
  in~\cite{GM}, or, simply {\it trisps} in~\cite{book}.} in a~standard
way: the vertices are indexed with numbers $1,\dots,d+1$, and the
order of the vertices in each boundary simplex is induced by that
global order. There is, in a~certain sense, a~unique way to identify
the boundary simplices with each other, if we want to identify as many
simplices as possible, while at the same time preserving the property
that the quotient has the induced $\da$-complex structure. This is
achieved by taking any two boundary simplices of the same dimension
and gluing them by the unique linear isomorphism which preserves the
order of the vertices. When $d=1$, we obtain a~circle. When $d=2$, we
obtain the so-called {\it Dunce hat}. This is a~classical
$\da$-complex, which is contractible but not collapsible. Its second
barycentric subdivision is a~simplicial complex that is of course also
contractible but not collapsible.

Graphically, we think about this boundary simplex identification as
{\it scrambling}.  Viewed as a~$\da$-complex, this {\it maximally
  scrambled} $d$-simplex will have a~single simplex in each
dimension. This paper grew out of the attempts by the author to better
understand and to generalize this well-known construction by relaxing
the {\it scrambling condition}. Our suggestion is that, guided by
a~certain combinatorial principle, we identify only some of the
boundary simplices of the same dimension, instead of gluing together
all of them.

Specifically, we start by putting labels on the vertices of the
original $d$-simplex. The label of vertex $i$ is denoted by $a_i$, for
all $i=1,\dots,d+1$. We think about these labels as {\it letters}, so
the ordered sequence of labels $(a_1,\dots,a_{d+1})$ gives a~{\it
  word}~$w$. Given any subsimplex $(i_1,\dots,i_t)$, where
$i_1<\dots<i_t$, of the original $d$-simplex, the corresponding
ordered label sequence $(a_{i_1},\dots,a_{i_t})$ gives a~{\it subword}
of~$w$. Using the same order-preserving linear isomorphism as above,
we now identify any two boundary simplices, which yield the same
subword. These simplices are necessarily of the same dimension, as
this is just the length of the subword minus one. It is easy to see
that the obtained space has the induced structure of a~$\da$-complex,
which we denote~$\da_w$.

The maximally scrambled case above corresponds to putting the same
label on all the vertices of $\da^d$, or, in other words, the word $w$
consists of a~single letter repeated $d+1$ times. At another extreme,
if all the labels are different, then no identifications take place at
all, and our space is just the original $d$-simplex itself. The whole
family $\{\da_w\}$ can then be seen as approximating between simplices
and higher-dimensional Dunce hats.

There is a rich and well-developed theory of combinatorially defined
simplicial complexes, see, e.g., \cite{book}. However, to our
knowledge very little work has been done in the category of
$\da$-complexes, see, e.g., \cite{Ko00,Ko02}, where combinatorially
defined $\da$-complexes, called there triangulated spaces, have been
applied to analyze spaces of polynomials with multiple
roots. Accordingly, we view the family $\{\da_w\}$ as a~source of
interesting combinatorially defined $\da$-complexes.

In this paper, we are able to completely determine the homotopy type
of the $\da$-complexes $\da_w$, see our main Theorem~\ref{thm:main}.
Somewhat surprisingly, all these spaces are either contractible or
homotopy equivalent to odd-dimensional spheres. The combinatorial rule
which reads off from the word $w$ the homotopy type of the
corresponding space~$\da_w$ is as follows. If
$w=a_1v_1a_1a_2v_2a_2\dots a_qv_qa_q$, where $a_1,\dots,a_q$ are (not
necessarily distinct) letters, and $v_1,\dots,v_q$ are (possibly
empty) words, such that $a_i$ does not occur in $v_i$, for all $1\leq
i\leq q$, then $\da_w$ is homotopy equivalent to $S^{2q-1}$, else
$\da_w$ is contractible.

As the added benefit, we are able to emulate the Dunce hat phenomenon,
and to obtain a~large family of both $\da$-complexes, as well as
simplicial complexes, which are contractible, but not collapsible.

\subsection{Preliminaries} $\,$

\nin
We use the notation $[n]=\{1,\dots,n\}$, for all natural numbers~$n$.
Let us recall the definition of the $\da$-complexes.

\begin{df}
The {\bf gluing data} which defines a~{\bf $\da$-complex} $X$ consists of
two parts:
\begin{itemize}
\item a family of sets $S_0,S_1,S_2,\dots$;
\item for each $0\leq m\leq n$ and each order-preserving injection
$f:[m+1]\hookrightarrow[n+1]$, we have a set map $B_f:S_n\ra S_m$. 
\end{itemize}
This data is subject to the following conditions:
\begin{enumerate}
\item[(1)] for any pair of composable order-preserving injections
$g:[k+1]\hra[m+1]$ and $f:[m+1]\hra[n+1]$, we have $B_{f\circ g}=B_g\circ B_f$;
\item[(2)] for any identity map $\id_n:[n+1]\hra[n+1]$, we have
  $B_{\id_n}=\id_{S_n}$.
\end{enumerate} 
\end{df}

The sets $S_n$ are the sets of {\it $n$-simplices} of $X$ and the maps
$B_f$ are the {\it boundary maps}. Assume now that for some $n\geq 1$,
we have picked $\sigma\in S_{n-1}$ and $\tau\in S_n$. For an
order-preserving injection $f:[n]\hra[n+1]$, set $\sgn f:=(-1)^k$,
where $k$ is the unique element in $[n+1]\sm\im f$. We now set
\[[\sigma:\tau]:=\sum_f \sgn f,\]
where the sum is taken over all order-preserving injections
$f:[n]\hra[n+1]$, such that $B_f(\tau)=\sigma$. We refer to
\cite[Section~2.1]{Ha02} and \cite[Section 2.3]{book} for further
general details on $\da$-complexes.

For the sake of being self-contained, we define the notions which we
need in this paper.  In the next two definitions, assume we are given
two $\da$-complexes $X$ and $\tilde X$ with respective gluing data
$\left(\{S_n\},\{B_f\}\right)$ and $\left(\{\tilde S_n\},\{\tilde
B_f\}\right)$.

\begin{df} \label{df:dais}
We say that $X$ is {\bf isomorphic} to $\tilde X$, if there exists
a~family of bijections $\alpha_n:S_n\ra\tilde S_n$, for all $n\geq 0$,
satisfying a~commuting relation
\[\tilde B_f\circ\alpha_n=\alpha_m\circ B_f:S_n\ra\tilde S_m,\]
for all order-preserving injections $f[m+1]\hookrightarrow[n+1]$.
\end{df}

Such a~family of bijections $\{\alpha_n\}$ is also called a~{\it
  $\da$-complex isomorphism} between the complexes $X$ and~$\ti X$.

\begin{df}
We define a new $\da$-complex, with gluing data $(\{T_n\},\{C_f\})$,
which we call the {\bf join} of $X$ and $\tilde X$, and denote by
$X*\tilde X$. To that end, we set
\[T_n:=\bigcup_{i+j+1=n}\{(\sigma,\tilde\sigma)\, |\,
\sigma\in S_i,\,\,\tilde\sigma\in\tilde S_j\}.\] Let now
$f:[m+1]\hookrightarrow[n+1]$ be an order-preserving injection, and
pick $\sigma\in S_i$, $\tilde\sigma\in\tilde S_j$, such that
$i+j+1=n$. The map $f$ can be represented by order-preserving
bijections $f:[\alpha+1]\hra[i+1]$, and $\tilde f:
[\tilde\alpha+1]\hra[j+1]$, where $\alpha:=|\im f\cap[i+1]|-1$, and
$\tilde\alpha:=|\im f\cap\{i+2,\dots,n+1\}]|-1$. We then set
$C_f((\sigma,\tilde\sigma)):=(B_f(\sigma),\tilde B_{\tilde
  f}(\tilde\sigma))$.
\end{df}


When $X$ and $Y$ are $\da$-complexes, we shall abuse notations and
denote by the same letters the corresponding CW complexes and
corresponding topological spaces, i.e., their geometric realizations.
There will be different ways by which we shall relate our objects to
each other. In connection with that we would like to remind the reader
about the following sequence of implications:
\[
\begin{array}{c}
\textrm{$X$ and $Y$ are isomorphic as $\da$-complexes}\\
\Downarrow\\
\textrm{$X$ and $Y$ are isomorphic as CW complexes}\\
\Downarrow\\
\textrm{$X$ and $Y$ are homeomorphic}\\
\Downarrow\\
\textrm{$X$ and $Y$ are homotopy equivalent.}
\end{array}
\]

There is a notion of {\it elementary collapses} for CW complexes, see
\cite[\S 4]{Co}. We do not need the full generality for
$\da$-complexes.  The following operation is very close to elementary
collapses of simplicial complexes and is sufficient for our purposes.

\begin{df} \label{df:dacoll}
Consider a~$\da$-complex $X$ given by the gluing data
$\left(\{S_n\},\{B_f\}\right)$. Assume we have $\tau\in S_n$, and $\sigma\in S_{n-1}$,
such that
\begin{enumerate}
\item[(1)] there exists a~unique order-preserving injection $f:[n]\hra[n+1]$, such that
  $B_f(\tau)=\sigma$;
\item[(2)] if we have $\tilde\tau\in S_n$, and an order-preserving injection $\tilde
  f:[n]\hra[n+1]$, such that $B_{\tilde f}(\tilde\tau)=\sigma$, then $f=\tilde f$ and
  $\tau=\tilde\tau$;
\item[(3)] the simplex $\tau$ is maximal in the following sense: there does not exist
  another simplex $\delta\in S_{n+1}$, such that $B_g(\delta)=\tau$, for some
  order-preserving injection $g:[n+1]\hra[n+2]$.
\end{enumerate}

Removing $\tau$ from $S_n$, $\sigma$ from $S_{n-1}$, and restricting the maps $B_f$
accordingly yields a~new $\da$-complex, which we shall call $X\sm\{\sigma,\tau\}$. We say
that it is obtained from $X$ by an {\bf elementary collapse}.
\end{df}
Sometimes, we call the pair $(\sigma,\tau)$ itself an elementary
collapse. When $X$ and $Y$ are $\da$-complexes, we have the following
sequence of strict implications:
\[
\begin{array}{c}
\textrm{there exists a sequence of elementary collapses reducing $X$ to $Y$}\\
\Downarrow\\
\textrm{there exists a~strong deformation retraction from $X$ onto $Y$}\\
\Downarrow\\
\textrm{$X$ and $Y$ are homotopy equivalent.}
\end{array}
\]

Finally, we note that all homology groups which we consider in this paper are taken with
integer coefficients.

\subsection{The scrambled simplices}\label{ssect:red} $\,$

\nin
Let us now describe the language in which we want to talk about the
scrambled simplices.

\begin{df}
Given any set $S$, we define a {\bf word} $w$ {\bf in alphabet $S$} to
be any finite ordered tuple $(a_1,\dots,a_n)$ of elements of $S$; we
allow repetitions in that tuple. The elements $a_1,\dots,a_n$ are
referred to as {\bf letters} of~$w$.  The number $n$ is called
the~{\bf length} of~$w$, which we denote by~$l(w)$. We set $\supp w:=
\{a_1,\dots,a_n\}$, and call it the {\bf support} set of $w$.
\end{df}
\noindent Note, that $|\supp w|\leq l(w)$, and in general it is possible to have
the strict inequality.
\vspace{1ex}

We shall write $w=a_1,\dots,a_n$, where for all $1\leq i\leq n$, we
have $a_i\in S$. Oftentimes we shall skip the commas and simply write
$w=a_1\dots a_n$, We shall use the power notation to denote
repetitions of letters, so $w=a^3$ means $w=aaa$, and $w=(a^2b)^2
b=a^2ba^2b^2=aabaabb$. For any $0\leq k\leq n$, we shall say that the
word $a_1\dots a_k$ is a~{\it prefix} of $w$; when additionally $k\leq
n-1$, we shall say that it is a~{\it proper prefix} of~$w$.

We say that $w=a_1^{\alpha_1}\dots a_t^{\alpha_t}$ is the {\it reduced
  form}\footnote{We find convenient to slightly abuse notations here
  and use the same letters $a_i$.} of the word $w$, if $a_i\neq
a_{i+1}$, for all $1\leq i \leq t-1$, and $\alpha_i\geq 1$, for all
$1\leq i\leq t$.  Obviously, every word has a unique reduced form,
and, when not stated otherwise, we shall assume that our words are
written in a~reduced form.

\begin{df}
Assume we are given a word $w=a_1\dots a_n$, and a subset
$I\subseteq[n]$, say $I=\{i_1,\dots,i_k\}$, where $i_1<\dots<i_k$.  We
set $w_I:=a_{i_1}\dots a_{i_k}$, and call it the {\bf $I$-subword}
of~$w$.
\end{df}

It is convenient to identify a subset $I\subseteq[n]$, such that
$|I|=k$, with an order-preserving injection $I:[k]\hra[n]$.

\begin{df}
Given a word $w=a_1\dots a_n$, $a_i\in S$, $n\geq 1$, the {\bf
  scrambled simplex} $\da_w$ is the $\da$-complex defined as follows.
\begin{itemize}
\item For all $l\geq 0$, we set $S_l$ to be the set of all
  $I$-subwords of $w$, such that $|I|=l+1$.
\item Assume we are given an order-preserving injection
  $f:[m+1]\hookrightarrow[l+1]$, where $0\leq m\leq l$. Take $w_I\in
  S_l$, where $|I|=l+1$. We have an order-preserving injection
  $I:[l+1]\hra[n]$, and can consider the composition $I\circ
  f:[m+1]\hra [n]$. We now set \[B_f(w_I):=w_{\im\,(I\circ f)}.\]
\end{itemize}
\end{df}

Clearly, the $\da$-complex $\da_w$ only depends on the underlying
partition of $[t+1]$ corresponding to the word $w$ in the following
sense: given words $w$ in alphabet $S$ and $w'$ in alphabet $S'$, and
a~renaming function $f:S\ra S'$, such that $f(w)=w'$, then $f$ induces
a~$\da$-complex isomorphism, see Definition~\ref{df:dais}, between
$\da_w$ and $\da_{w'}$.

Intuitively, the gluing data of $\da_w$ simply records what happens
when we delete letters. As mentioned above, there is an alternative
description of $\da_w$ as a quotient complex of the $d$-simplex
$\da^d$, where $d=l(w)-1$. In this description, we start with
a~$d$-simplex $\da$. Its boundary simplices are indexed by the subsets
of $[d+1]$, so let $\da^I$ denote the boundary simplex corresponding
to $I\subseteq[d+1]$. Now, if $w_I=w_J$, then we identify $\da^I$ with
$\da^J$ using a~linear isomorphism which preserves the order of the
vertices. In particular, if $w=a_1\dots a_{d+1}$, for $a_i\neq a_j$,
for all $i,j$, then $\da_w$ is just a~$d$-simplex. Accordingly, we
think of $\da_w$ as a~$d$-simplex whose boundary has been scrambled
in a~certain pattern, given by the word~$w$.

It is immediate that $\dim\da_w=d$. Furthermore, we have a~cellular
isomorphism $\da_w\simeq\da_{\bar w}$, where $\bar w$ is the word $w$
written backwards. This isomorphism does not have to be
a~$\da$-complex isomorphism, but it certainly induces
a~homeomorphism. For example, we have
$\da_{abb}\cong\da_{bba}\cong\da_{aab}$.

\section{Examples and first properties}

\subsection{The $\da$-complexes of subwords of a word of length at most 3} $\,$

\nin
We shall now describe, up to isomorphism, the cell complexes $\da_w$,
when $l(w)\leq 3$.  If $l(w)=1$, then we only need to consider $w=a$,
and clearly $\da_a$ is just a point.

When $l(w)=2$, we have two cases: $w=ab$ and $w=aa$. We see that
$\da_{ab}$ is a~$1$-simplex, and $\da_{aa}$ is the CW complex with one
$0$-cell and one $1$-cell, which is homeomorphic to $S^1$.

\begin{figure}[hbt]

  \input{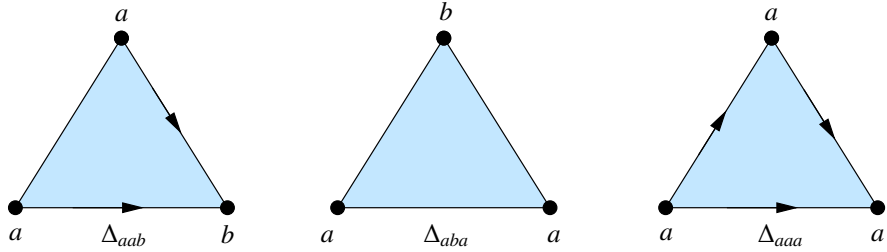}  

\caption{The $\da$-complexes of subwords of a word of length at most
  three. The arrows on the edges indicate which edges are glued
  together and in which direction the gluing is done.}
\label{fig:1}
\end{figure}

When $l(w)=3$, we have the cases: $w=abc$, $w=a^2b$, $w=aba$, and
$w=a^3$. Accordingly, $\da_{abc}$ is a $2$-simplex, $\da_{a^2b}$ is
homeomorphic to a~disc, so is contractible, and $\da_{aba}$ is
homeomorphic to the topological space obtained from the disc by
identifying two of its boundary points, which is homotopy equivalent
to $S^1$. Finally, $\da_{a^3}$ is the classical Dunce hat, see, e.g.,
\cite{Ze64} for further details. It is well-known that this space is
contractible as well. All the three nontrivial cases are shown on
Figure~\ref{fig:1}.

\subsection{Concatenation of words from disjoint alphabets} $\,$

\nin
Before we proceed with computing further examples, we make the
following simple, but useful proposition, whose formal verification is
left to the reader.
\begin{prop}
Assume $w$ is a~concatenation of two words $w=w_1\circ w_2$, such that
$\supp{w_1}$ and $\supp{w_2}$ are disjoint. Then we have
a~$\da$-complex isomorphism.
\end{prop}

\begin{equation} \label{eq:concat}
\Delta_w\simeq\Delta_{w_1}*\Delta_{w_2}.
\end{equation}

In a situation like this, we shall say that $w$ is {\it decomposable},
and else we say that $w$ is {\it indecomposable}. The
Tables~\ref{table:ind1} and~\ref{table:ind2} list, up to renaming and
up to reversing the order of the letters, all indecomposable words of
length at most~$5$.

{\renewcommand{\arraystretch}{1.3}
\begin{table}[hbt] 
\[\begin{array}{l|l}
l(w) & \text{indecomposable words}\\ \hline
1 & a \\ \hline 
2 & a^2 \\ \hline
3 & a^3, aba \\ \hline
4 & a^4, aba^2, abab, ab^2a, abca \\ \hline
\end{array}\]
\caption{The $9$ indecomposable words of length at most~$4$.}
\label{table:ind1}
\end{table}
\begin{table}[hbt] 
\[\begin{array}{l|l}
|\supp{w}| & \text{indecomposable words}\\ \hline
1 & a^5 \\ \hline
2 & aba^3, a^2ba^2, abab^2, ab^2ab, ab^3a, ab^2a^2, ababa \\ \hline
3 & abca^2, abaca, abacb, ab^2ca \\ \hline
4 & abcda \\ \hline
\end{array}\]
\caption{The $13$ indecomposable words of length~$5$.}
\label{table:ind2}
\end{table}
}

Clearly, when interested in determining the topology
of the complexes $\da_w$ it is fully sufficient to restrict ourselves
to considering the indecomposable words only.

\subsection{The higher-dimensional Dunce hats} $\,$

\nin
The indecomposable words $w=a^n$ correspond to an interesting family
of spaces. As mentioned above $\da_{a^3}$ is a classical Dunce hat.
It is then easy to understand $\da_{a^4}$. Indeed, the $2$-skeleton of
this CW complex is homeomorphic to a Dunce hat, hence it is
contractible. Therefore, $\da_{a^4}$ is obtained by attaching
a~$3$-cell to a~contractible space.  Contracting the $2$-skeleton to
a~point we see that $\da_{a^4}$ is homotopy equivalent to $S^3$.

In general, the spaces $\da_{a^n}$ were defined by Andersen,
Marjanovi{\'c}, and Schori, see~\cite{AMS}, using the symmetric
products of intervals, see also Borsuk and Ulam, \cite{BU}.
When $n$ is odd, the spaces $\da_{a^n}$ are called
{\it higher-dimensional Dunce hats}. 
The topology of $\da_{a^n}$ was completely determined in~\cite{AMS}.

\begin{prop}\label{prop:an}{\rm (\cite[Theorem~2.3]{AMS}).}
When $n$ is odd, the spaces $\da_{a^n}$ are contractible, and when $n$
is even, the space $\da_{a^n}$ is homotopy equivalent to a sphere of
dimension $n-1$. 
\end{prop}
The case when $n$ is odd is the interesting one and was proved using
Homotopy Addition Theorem of Hu, \cite{Hu}. The case when $n$ is even
is a~simple corollary. Since $\da_{a^n}$ is always obtained from
$\da_{a^{n-1}}$ by attaching a single $(n-1)$-cell, our above argument
for $n=4$ works in general. The space $\da_{a^{n-1}}$ is contractible,
and its inclusion into $\da_{a^n}$ is a~cofibration, see,
e.g.,~\cite{book}, so we can just shrink it to a~point without
changing the homotopy type, and end up with an~$(n-1)$-sphere.

Alternatively, it is easy to prove Proposition~\ref{prop:an} using
a~version of Whitehead's theorem, as is done in this paper.

\section{Formulas for the reduced Euler characteristics}

\subsection{The arrow terminology and the recursive formula} $\,$

\nin In this section, we compute the Euler characteristic of the complexes
$\da_w$ for all words~$w$. This simple derivation is the first step in the
general understanding of the homotopy type of the spaces $\da_w$.  We also
design the combinatorial language, which we will need later to formulate our
main result.

\begin{df}
For any word $w$, let $\ce(w)$ denote the reduced Euler
characteristics of the $\da$-complex~$\da_w$.  
\end{df}

Note that $\ce(w)$ counts the subwords of $w$ with weights, with each
word giving the contribution $(-1)^{l(w)+1}$. For example, we have
$\ce(a^{2t})=-1$, $\ce(a^{2t+1})=0$, $\ce(aba)=-1$, etc. Also, for the
empty word we have $\ce(\emptyset)=-1$. We shall prove the following,
somewhat surprising result:

\begin{center}
{\it for all words $w$, we either have $\ce(w)=0$, or $\ce(w)=-1$.}
\end{center}

\noindent
The more precise statement is given in Theorem~\ref{thm:ce}.  

Before we can prove the result on the Euler characteristics, we need to
introduce some new terminology.

\begin{df}
Assume $w$ is a word, and pick a~letter $a\in\supp w$.  The new word, which we
denote $w\dar a$, is obtained from $w$ by finding the leftmost occurrence of $a$
in $w$ and deleting everything to the left of it, including $a$ itself.
\end{df}

For example, $a^n\dar a=a^{n-1}$, $aba\dar a=ba$, $aba\dar b=a$, and if $w=va$,
such that $a\notin\supp v$, then $w\dar a=\emptyset$. Let us furthermore
introduce the following short hand notation: for any word $v=a_1\dots a_t$,
where $a_i$'s are letters, we set
\[w\dar v:=(\dots(w\dar a_1)\dar a_2)\dar\dots)\dar a_t.\]

\begin{prop}\label{prop:wx}
Assume $w$ is a word, and $a$ is a letter in its support. There is
a~1-1-correspondence between the subwords of $w\dar a$ and those subwords of $w$
which begin with the letter~$a$. This correspondence is given by adding $a$ as
the first letter to a~subword of $w\dar a$.
\end{prop}
\pr Let $A$ denote the set of all subwords of $w\dar a$, and let $B$ denote the
set of all subwords of $w$ which begin with the letter~$a$.  Let $\varphi:A\ra
B$ be the map which adds $a$ as the first letter, and let $\psi:B\ra A$ be the
map which deletes the first letter (which, by definition of $B$, must be~$a$).

It is obvious that $\varphi$ is well-defined, since $a$ can always be added on
the left to any subword of $w\dar a$. The $\psi$ is well-defined is equally
obvious, since when $av$ is a~subword of $w$, $v$ will be the subword of $w\dar
a$, essentially by the definition of $w\dar a$. The two maps are also inverses
of each other, so we are done with the proof.  \qed

\vspace{5pt}

\nin The next proposition allows us to use recursion to calculate the
function~$\ce$.

\begin{prop} \label{prop:rce}
Assume $w$ is an arbitrary word, then we have the following recursive
formula
\begin{equation}\label{eq:rce}
\ce(w)=-\sum_{x\in\supp{w}}\ce(w\dar x)-1.
\end{equation}
\end{prop}
\pr Recall that the left hand side of \eqref{eq:rce} counts the
subwords of $w$, with each word giving the contribution
$(-1)^{l(w)+1}$. The empty word gives the term $-1$, and all other
words start with some letter from $\supp w$. Hence we have
\begin{equation}\label{eq:rcea}
\ce(w)=-1+\sum_{x\in\supp{w}}E_x,
\end{equation}
where for each $x\in\supp w$, the term $E_x$ denotes the total contribution of
the subwords of $w$ which start with~$x$. By Proposition~\ref{prop:wx}, there is
a~bijection between the set of such words and the subwords of $w\dar x$. This
bijection changes the length of the word by~$1$, hence it changes the sign of
the contribution. We can therefore conclude that for all $x\in\supp w$, we have
the equality $E_x=-\ce(w\dar x)$. Substituting this into~\eqref{eq:rcea} we
obtain the identity~\eqref{eq:rce}.  \qed

\vspace{5pt}

In particular, if $\supp{w}=\{a\}$, then for all $n\geq 1$, we get
$\ce(a^n)=-\ce(a^{n-1})-1$, which we can rewrite as
$\ce(a^{n-1})+\ce(a^n)=-1$. This is consistent with the direct observation that
\begin{equation} \label{eq:cir1}
\ce(a^n)=\begin{cases} 
-1, & \text{ if }n\text{ is even};\\
0, & \text{ if }n\text{ is odd}.
\end{cases}
\end{equation}

If $\supp{w}=\{a,b\}$, then \eqref{eq:rce} says
\begin{equation}\label{eq:rce2}
\ce(w)=-\ce(w\dar a)-\ce(w\dar b)-1.
\end{equation}
Let $\alt(t)$ denote the word of length $t$ consisting of two alternating
letters, say $\alt(4)=abab$. As another example, we calculate
$f_t:=\ce(\alt(t))$. Clearly, we have $f_0=-1$ and $f_1=f_2=0$. Equation
\eqref{eq:rce2} tells us that $f_t=-f_{t-1}-f{t-2}-1$, for all $t\geq
2$. Reindexing and moving terms yields the following identity:
\begin{equation}
\label{eq:fn}
f_n+f_{n+1}+f_{n+2}=-1,\text{ for all }n\geq 0.
\end{equation}
Comparing \eqref{eq:fn} for two consecutive values of $n$, we conclude that
$f_n$ is periodic with period $3$, in other words that $f_n=f_{n+3}$, for all
$n\geq 0$. We conclude that

\begin{equation} \label{eq:cir2}
\ce(\alt(n))=\begin{cases} 
-1, & \text{ if }n\text{ is divisible by }3;\\
0, & \text{ otherwise}.
\end{cases}
\end{equation}
Later we shall see that $\da_{\alt(n)}$ is actually contractible
unless $n$ is divisible by $3$. If $n=3k$, then we will show that
$\da_{\alt(n)}\simeq S^{2k-1}$.

In fact, the equations~\eqref{eq:cir1} and~\eqref{eq:cir2} can easily be
generalized to words, where a certain set of $t$ letters repeats in a~circular
manner: $w=a_1\dots a_ta_1\dots a_t\dots$. For such a word we have $\ce(w)=-1$
if the length of $w$ is divisible by $t+1$, and $\ce(w)=0$ otherwise.

\subsection{Elimination of circular words and the main Euler characteristics 
theorem} $\,$

\nin The following proposition provides the crucial step in our
computation.

\begin{prop}\label{prop:red2} 
Assume $w=auav$, where $a$ is a letter, and $u$ and $v$ are words,
such that $a\notin\supp{u}$. Then we have the equality
$\ce(w)=\ce(v)$.
\end{prop}

\begin{rem}
We would like to point out that it is allowed for words $u$ and $v$ in
Proposition~\ref{prop:red2} to be empty. When $v$ is empty, we recover
a simple corollary of Proposition~\ref{prop:ava}. When $u$ is empty,
we obtain the equality $\ce(a^2 v)=\ce(v)$, for an arbitrary word~$v$.
\end{rem}

\nin {\bf Proof of Proposition~\ref{prop:red2}.}

\nin
By the recursive formula \eqref{eq:rce}, we have 
\begin{equation}\label{eq:rd1}
\ce(w)=-\sum_{\substack{x\in\supp{w}\\ x\neq a}}\ce(w\dar x)-\ce(w\dar a)-1.
\end{equation} 
Note that $w\dar a=uav$. Applying~\eqref{eq:rce} to that word, we obtain
\begin{equation}\label{eq:rd2}
\ce(uav)=-\sum_{\substack{x\in\supp{uav}\\x\neq a}}\ce(uav\dar x)-\ce(uav\dar a)-1.
\end{equation}
Let us now substitute~\eqref{eq:rd2} into~\eqref{eq:rd1}.
We get
\begin{multline} \label{eq:mult2}
\ce(w)=-\sum_{\substack{x\in\supp{w}\\x\neq a}}\ce(w\dar x)+
\left(\sum_{\substack{x\in\supp{uav}\\x\neq a}}\ce(uav\dar x)+\ce(uav\dar a)+1\right)-1=\\
=\ce(uav\dar a)+\sum_{\substack{x\in\supp{w}\\x\neq a}}\left(\ce(uav\dar x)-\ce(w\dar x)\right)=\ce(v),
\end{multline} 
where the penultimate equality is obtained by using the fact that
$\supp{w}=\supp{uav}$, and the last equality follows from the equalities $uav\dar a=v$,
and $uav\dar x=w\dar x$, whenever $x\in\supp w$, $x\neq a$.
\qed

\begin{df}
Let $w$ be an arbitrary word. We introduce the following terminology.
\begin{itemize}
\item We call $w$ {\bf circular}, if it is of the form $w=ava$, where
  $a$ is a letter, and $v$ is a~word, which is possibly empty, such
  that $a\notin\supp{v}$.
\item We call $w$ {\bf spherical} if it can be represented as
  a~concatenation of circular words. 
\item We call $w$ {\bf conical} if it is of the form $w=av$, where $a$
  is a letter, and $v$ is a~word, which is possibly empty, such that
  $a\notin\supp{v}$.
\end{itemize}
\end{df}

Note, that we consider the empty word to be spherical, and we view it
as a~concatenation of the empty set of circular words. We do not
consider the empty word to be either circular, or conical.

Clearly, when $w$ is conical, the $\da$-complex $\da_w$ is a cone with
apex $a$; in particular, $\ce(w)=0$. Furthermore, if $w=ava$ is
circular, then Proposition~\ref{prop:ava} tells us that $\da_w$ is
homotopy equivalent to a~circle, and so $\ce(w)=-1$. This explains our
terminology.

\begin{prop}\label{prop:dec} The following is true for all words.
\begin{enumerate}
\item[(1)] A decomposition of a spherical word into circular one is unique.
\item[(2)] If a~word $w$ is not spherical, then there is a unique
  decomposition $w=uv$, where $u$ is a~spherical word, which is
  possibly empty, and $v$ is a~conical word.
\end{enumerate}
\end{prop}
\pr To see (1), let $a$ be the first letter of a~spherical word $w$.
The first circular word in the decomposition of $w$ must be the
subword between the two leftmost occurrences of $a$, including the
letter $a$ on both ends. Proceeding left to right we see that the
entire decomposition is unique.

To see (2) proceed as in the argument above. Clearly, unless the word
is spherical, we will end up with the decomposition $w=w_1\dots
w_ka\tilde v$, where the words $w_1,\dots, w_k$ are circular, and $a$
is a letter, such that $a\notin\supp{\tilde v}$. Setting $u:=w_1\dots
w_k$, and $v:=a\tilde v$, we obtain the desired decomposition.
\qed

\vspace{5pt}

As the first example, consider the words $w=a^n$. When $n=2$, such a
word is circular. Hence, when $n$ is even, such a word is
spherical. When $n$ is odd, we can decompose $w=a^{n-1}\cdot a$, where
$a^{n-1}$ is spherical and $a$ is conical.

Another example is provided by the alternating words $w=\alt(n)$.  The
decomposition $(ab)^3=aba\cdot bab$ shows that $(ab)^3$ is spherical,
and hence $(ab)^{3t}$ is spherical for any~$t$. When $n=3t+1$, we get
the decomposition $w=(ab)^{3t}\cdot ab$, where $ab$ is conical. When
$n=3t+2$, we get the decomposition $w=(ab)^{3t}\cdot aba \cdot b$,
where $b$ is conical, and the word before it is spherical.

\begin{thm}\label{thm:ce}
For an arbitrary word $w$, we have
\begin{equation} 
\ce(w)=\begin{cases}
-1, & \textrm{ if } w\textrm{ is spherical;} \\
 0, & \textrm{ otherwise.}
\end{cases}
\end{equation}
\end{thm}
\pr If $w$ is spherical, then it follows from repeated application of
Proposition~\ref{prop:red2} that $\ce(w)=\ce(\emptyset)=-1$. 
If $w$ is not spherical, by Proposition~\ref{prop:dec}(2), we can
write $w=uv$, where $u$ is spherical and $v$ is conical. Repeatedly
applying Proposition~\ref{prop:red2}, we obtain $\ce(w)=\ce(v)=0$.
\qed

\section{Fundamental group of $\da_w$}

\nin We start with a simple observation.

\begin{prop}\label{prop:ava}
When the word $w$ is circular, the $\da$-complex $\da_w$
is homotopy equivalent to $S^1$.
\end{prop}
\pr Assume $w=ava$, where $a$ is a letter, and $v$ is a~word, such
that $a\notin\supp{v}$.  Let $\tilde a$ be a letter, such that $\tilde
a\notin\supp{w}$, and set $\tilde w:=av\ti a$. By the concatenation
formula~\eqref{eq:concat}, the cell complex $\da_{\ti w}$ is obtained
from $\da_v$ by coning twice, with apexes $a$ and $\tilde a$. In
particular, the space $\da_{\ti w}$ is contractible. On the other
hand, the $\da$-complex $\da_w$ is obtained from $\da_{\ti w}$ by
identifying two of its vertices, namely the ones labeled by $a$ and by
$\ti a$. The proposition follows now from an easy general fact, that
gluing together two vertices in a~connected cell complex $K$ results
in a~space which is homotopy equivalent to a wedge of $K$
with~$S^1$. \qed

\vspace{5pt}

Proposition~\ref{prop:ava} covers many cases of the words from
Tables~\ref{table:ind1} and~\ref{table:ind2}. When $l(w)=4$, it covers
words $ab^2a$ and $abca$, and when $l(w)=5$, we get the words $ab^3a$,
$ab^2ca$, and $abcda$.

It turns out, that not only the circular words are the only ones for
which $\da_w$ is homotopy equivalent to $S^1$, but, in fact, these are
the only words, for which $\da_w$ is not simply connected.

\begin{thm}\label{thm:fund}
Let $w$ be an arbitrary word.  The fundamental group of the
$\da$-complex $\da_w$ is given by 
\[ \pi_1(\da_w)\simeq
\begin{cases}
\zz, & \text{  if }w\text{ is circular;} \\
0, & \text{ otherwise.}
\end{cases}
\]
\end{thm}
\pr By Proposition~\ref{prop:ava} we know that when $w$ is circular, the $\da$-complex
$\da_w$ is actually homotopy equivalent to a circle, so obviously,
$\pi_1(\da_w)\simeq\zz$.

Assume now that the word $w=a_1\dots a_t$ is not circular. Take $x:=a_1$ to be the base
point, and consider $\pi_1(\da_w,x)$. To start with, if $x\notin\{a_2,\dots,a_t\}$, then
$\da_w$ is a cone with apex $x$, so the fundamental group is trivial. Assume therefore,
that $x\in\{a_2,\dots,a_t\}$. Let $m$ denote the minimal index such that $m\geq 2$ and
$a_m=x$. Since the word $w$ is not circular, we also have $m\leq t-1$.

The subword $a_1a_m=xx$ indexes an edge of $\da_w$, which is a~loop based at $x$. If we
fix a~choice of orientation on that loop, we obtain a~representation of an element
$h_x\in\pi_1(\da_w,x)$. That element is easy to understand. Namely, consider the subword
$a_1 a_m a_t= xxa_t$. Independently of the fact, whether $a_t=x$, the boundary of the
corresponding $2$-simplex, viewed as a~loop in $\da_w$ tells us that $h_x=0$.

It is now a general fact about $\da$-complexes, that any element of the fundamental group
$g\in\pi_1(\da_w,x)$ can be represented as a~sequence of edges $xb_1,\dots,b_kx$, such
that $k\geq 0$, and $b_i\in\supp w$, for all $1\leq i\leq k$. Call this an {\it edge
  representation}.

Assume that the fundamental group $\pi_1(\da_w,x)$ is non-trivial.
Consider all nontrivial elements of $\pi_1(\da_w,x)$, and all their
edge representations. Pick among all these representations one which
minimizes $k$, say it consists of $l+1$ edges, and let $g$ denote the
represented element. Since the loop $h_x$ represents a trivial
element, we must have $l\geq 1$. Furthermore, since $l$ is minimal
possible, we have $x\notin\{b_1,\dots,b_l\}$.

Assume first $l=1$. We have $b_1\neq x$. Either $a_1b_1a_m=xb_1x$ or $a_1a_mb_1=xxb_1$ is
a~subword of~$w$.  The corresponding $2$-simplex provides a~path homotopy between the
concatenation of $xb_1$ and $b_1x$ and one of the orientations of the loop $xx$. Since the
latter represents a~trivial element of the fundamental group, we conclude that $g=0$,
yielding a~contradiction.

Assume now, that $l\geq 2$. We have $x\notin\{b_1,\dots,b_l\}$.  Either $xb_1b_2$ or
$xb_2b_1$ is a~subword of $w$, and therefore indexes a~$2$-simplex of~$\da_w$. This means
that the concatenation of edges $x b_1$ and $b_1 b_2$ is path-homotopic to the edge $x
b_2$, contradicting the assumption that $l$ is smallest possible.

We conclude that the group $\pi_1(\da_w,x)$ is trivial.  \qed

\section{The main theorem}

\subsection{Orders on exponential presentations of subwords} $\,$

\nin Assume $a_1^{\alpha_1}\dots a_t^{\alpha_t}$ is the reduced form
of a~word~$w$, see Subsection~\ref{ssect:red}. All of the subwords
of~$w$ can be written as $v=a_1^{\beta_1}\dots a_t^{\beta_t}$, where
$0\leq\beta_i\leq\alpha_i$, for all $i=1,\dots,t$; note that we are
forced to allow $\beta_i=0$. This is not necessarily a~reduced form of
the word $v$ and it is clearly not unique. For example, when
$w=aba=a^1 b^1 a^1$, we have two presentations for the subword $v=a$;
namely, $v=a^1 b^0 a^0$ and $v=a^0 b^0 a^1$.

Once the word $w$ fixed, it is sufficient to simply write the tuples
of the exponents to (non-uniquely) record the subwords.  So, in the
previous example, we could use the tuple $(1,1,1)$ to encode $w$
itself, while the two presentations of $v$ would be denoted by the
tuples $(1,0,0)$ and $(0,0,1)$.  We shall call such a tuple the {\it
  exponential presentation} of the word $v$ as a~subword of~$w$, and
we let $\cf_w(v)$ denote the set of all these exponential
presentations. So, in the above example, we have
$\cf_{aba}(a)=\{(1,0,0), (0,0,1)\}$.

The set of $t$-tuples of numbers may be equipped with various standard
orders. Here we will need the {\it domination order}, which is
a~partial order and is denoted by $>$, and two total orders: the~{\it
  lexicographic order}, denoted by $\succl$ and the~{\it
  colexicographic order}, denoted by~$\succc$. Let us recall what
these orders are.
\begin{itemize}
\item In the domination (partial) order, we say that
  $(x_1,\dots,x_t)\geq(y_1,\dots,y_t)$ if and only if $x_i\geq y_i$,
  for all $i=1,\dots,t$.
\item In the lexicographic order, we order the tuples as words in
  a~dictionary. In other words,
  $(x_1,\dots,x_t)>_{lex}(y_1,\dots,y_t)$ if and only if there exists
  $1\leq k\leq t$, such that $x_1=y_1,\dots,x_{k-1}=y_{k-1}$, and
  $x_k>y_k$.
\item The colexicographic order is essentially the same as the
  lexicographic one, except we read the $t$-tuples from right to left
  instead.  That is, $(x_1,\dots,x_t)>_{colex}(y_1,\dots,y_t)$ if and
  only if there exists $1\leq k\leq t$, such that
  $x_t=y_t,\dots,x_{k+1}=y_{k+1}$, and $x_k>y_k$.
\end{itemize}
Clearly, all of these orders are inherited by the set $\cf_w(v)$, for
arbitrary $w$ and~$v$.

\subsection{Shifted presentations} $\,$

\nin The next definition shall help us to standardize the ways we deal
with subwords of a~word.

\begin{df}
Assume, we are given a word $w$, whose reduced form is
$a_1^{\alpha_1}\dots a_t^{\alpha_t}$, and we are given
a~subword~$v$. We say that $(\beta_1,\dots, \beta_t)\in\cf_w(v)$ is
{\bf left-shifted} if it is maximal in $\cf_w(v)$ with respect to the
lexicographic order. We say that it is {\bf right-shifted} if it is
maximal in $\cf_w(v)$ with respect to the colexicographic order.
\end{df}

Obviously, for an arbitrary word $w$ and a subword $v$, both the
left-shifted and the right-shifted presentations of $v$ exist and are
unique.

\begin{df} 
Assume again that a~word $w$ is given by its reduced form
$a_1^{\alpha_1}\dots a_t^{\alpha_t}$, and let $p$ be an arbitrary
index, $1\leq p\leq t$.  Set $w':=a_1^{\alpha_1}\dots a_p^{\alpha_p}$,
and $w'':=a_p^{\alpha_p}\dots a_t^{\alpha_t}$.  We say that the tuple
$(\beta_1,\dots, \beta_t)\in\cf_w(v)$ is {\bf $p$-shifted} if
$(\beta_1,\dots, \beta_p)\in\cf_{w'}(v')$ is {\bf left-shifted} and
$(\beta_p,\dots, \beta_t)\in\cf_{w''}(v'')$ is {\bf right-shifted},
where $v'=a_1^{\beta_1}\dots a_p^{\beta_p}$, and
$v''=a_p^{\beta_p}\dots a_t^{\beta_t}$.
\end{df}

In the special cases $p=1$ and $p=t$, we recover the notions of being
right-shifted and being left-shifted, respectively. In general, in
contrast to the left- and right-shifted cases, we can only conclude
that a~$p$-shifted presentation exists. For example, a~$p$-shifted
presentation can easily by constructed by starting with any tuple
$(\beta_1,\dots,\beta_t)$ in $\cf_w(v)$, then first left-shifting the
tuple $(\beta_1,\dots,\beta_p)$, and then right-shifting the tuple
$(\beta_p,\dots,\beta_t)$. However, the $p$-shifted tuple is by no
means unique. In the example above, where $w=aba$ and $v=a$, both
$(1,0,0)$ and $(0,0,1)$ are $2$-shifted.  We do however have the
following proposition.

\begin{prop}\label{prop:pshift}
Assume $w$ is a~word, given by its reduced form $a_1^{\alpha_1}\dots
a_t^{\alpha_t}$, and let $v$ be a~subword of~$w$. Assume we have an
index $1\leq p\leq t$, and a~$p$-shifted exponential presentation of
$v$, $\beta=(\beta_1,\dots,\beta_t)$, such that $\beta_p\geq 1$. Then
the $p$-shifted exponential presentation of $v$ is unique.
\end{prop}
\pr Assume that we have another $p$-shifted exponential presentation of $v$, say
$\ti\beta=(\ti\beta_1,\dots,\ti\beta_t)$. Let $(\beta_{i_1},\dots,\beta_{i_k})$ be the
tuple of all $\beta_{i_1},\dots,\beta_{i_k}\neq 0$, such that $1\leq i_1<\dots<i_k\leq
p-1$, and $a_{i_1},\dots,a_{i_k}\neq a_p$. Symmetrically on the right, let
$(\beta_{j_1},\dots,\beta_{j_m})$ be the tuple of all $\beta_{j_1},\dots,\beta_{j_m}\neq
0$, such that $p+1\leq j_1<\dots<j_m\leq t$, and $a_{j_1},\dots,a_{j_m}\neq a_p$.  Let the
tuples $(\ti\beta_{\ti i_1},\dots,\ti\beta_{\ti i_{\ti k}})$ and $(\ti\beta_{\ti
  j_1},\dots,\ti\beta_{\ti j_{\ti m}})$ be defined the same way for $\ti\beta$.

Assume $a_{i_1}^{\beta_{i_1}}\dots a_{i_k}^{\beta_{i_k}}\neq a_{\ti i_1}^{\ti\beta_{\ti
    i_1}}\dots a_{\ti i_{\ti k}}^{\ti\beta_{\ti i_{\ti k}}}$.  Flipping the word, if
necessary, we can assume, without loss of generality, that $a_{i_1}^{\beta_{i_1}}\dots
a_{i_k}^{\beta_{i_k}}$ is a~{\it proper} prefix of $a_{\ti i_1}^{\ti\beta_{\ti i_1}}\dots
a_{\ti i_{\ti k}}^{\ti\beta_{\ti i_{\ti k}}}$. But that would mean that the whole subword
$a_1^{\beta_1}\dots a_p^{\beta_p}$ would have an exponential presentation within
$(\alpha_1,\dots,\alpha_{i_k})$. This contradicts to our assumption that the exponential
presentation $(\beta_1,\dots,\beta_p)$ is left-shifted.

We conclude that $(\beta_{i_1},\dots,\beta_{i_k})=(\ti\beta_{\ti i_1},\dots,\ti\beta_{\ti
  i_{\ti k}})$ and $(\beta_{j_1},\dots,\beta_{j_m})=(\ti\beta_{\ti
  j_1},\dots,\ti\beta_{\ti j_{\ti m}})$. Therefore, the presentations $\beta$ and
$\ti\beta$ may only differ on the exponents of those $a_i$, for $i_k<i<j_1$, for which
$a_i=a_p$. However, since $\beta_p>0$, we must have $\beta_i=\alpha_i$, whenever
$i_k<i<j_1$, $i\neq p$, and $a_i=a_p$. This obviously determines the rest of the
exponential presentation uniquely, and we conclude that $\beta=\ti\beta$.  \qed




\subsection{Whitehead's and Hurewicz' Theorems} $\,$

\nin
In this short subsection we list two classical results which we use
for our computations. Both are obtained by combining versions
of Whitehead's and Hurewicz' theorems.

\begin{thm} \label{thm:wh1}
A simply connected CW complex $X$ whose homology groups $\widetilde
H_i(X;\zz)$ are trivial is contractible.
\end{thm}

\begin{thm} \label{thm:wh2}
Let $\varphi:X\ra Y$ be a~map between simply connected CW complexes
that induces isomorphism maps $\varphi_*:H_n(X;\zz)\ra H_n(Y;\zz)$,
for all~$n$. Then, the map $\varphi$ is a~homotopy equivalence.  If,
furthermore, $\varphi$ is an inclusion map, then there exists a~strong
deformation retraction from $Y$ to~$X$.
\end{thm}

\nin We refer to various sources, such as \cite{W}, \cite[Corollary
  6.32 and Proposition 6.34]{book}, \cite[Theorem 4.4.5 and Corollary
  4.4.33]{Ha02}, and \cite[Theorem 7.6.25]{Sp}.

\subsection{Algebraic Morse Theory} $\,$

\nin In this subsection we present a short extract from the algebraic
Morse theory, which is sufficient for our purposes. We include
sketches of proofs to stay self-contained. The reader is advised to
consult~\cite[Section 11.3]{book}, \cite{dmt}, and the references
therein for a~more complete picture. For what follows, we recall that
when $P$ is a~partially ordered set, and $x\in P$, we set $P_{\geq x}
:=\{y\in P\,\,|\,y\geq x\}$, and $P_{>x}:=\{y\in P\,\,|\,y>x\}$.

\begin{df}\label{df:collapsing}
Assume we are given a~finite $\da$-complex $X$. Let $\cp(X)$ denote its face
poset. A~sequence of pairs of simplices $((\sigma_1,\tau_1),\dots,(\sigma_n,\tau_n))$ of
$X$ is called a~{\bf collapsing order} if the following conditions are satisfied:
\begin{enumerate}
\item[(1)] for all $1\leq i\leq n$, we have $\dim\sigma_i=\dim\tau_i-1$;
\item[(2)] $[\sigma_i:\tau_i]=\pm 1$;
\item[(3)] $\cp(X)_{\geq\sigma_i}\subseteq\{\sigma_1,\dots,\sigma_i,
\tau_1,\dots,\tau_i\}.$
\end{enumerate}
\end{df}

\begin{lm} \label{lm:dmt}
Assume we have a $\da$-complex $X$ of dimension $d+1$, and a pair of simplices
$(\sigma,\tau)$, such that $\dim\tau=d+1$, $\dim\sigma=d$, $[\sigma:\tau]=\pm 1$, and
$\cp(X)_{>\sigma}=\{\tau\}$. Let $\wti X$ be the $\da$-complex obtained from $X$ by
removing $\sigma$ and $\tau$.

Then, the inclusion map $\iota:\wti X\ra X$ induces isomorphism on homology groups
with integer coefficients. If, in addition, the spaces $X$ and $\wti X$ are simply
connected, then there is a strong deformation retraction from $X$ to $\wti X$.
\end{lm}

\pr A direct analysis of the chain complex $C_*(X,\wti X)$ shows that the homology
groups $H_n(X,\wti X)$ are trivial for all $n$. The long exact sequence of the pair
$(X,\wti X)$ then implies that $\iota$ induces isomorphism on homology groups.
Furthermore, the statement about the strong deformation retraction is a~direct
corollary of Theorem~\ref{thm:wh2}.  \qed

\begin{thm}\label{thm:dmt}
Assume we have a finite $\da$-complex $X$, and a~collapsing sequence
$(\sigma_1,\tau_1),\dots,(\sigma_q,\tau_q)$. Let $\wti X$ be the
$\da$-complex obtained from $X$ by removing the set of simplices
$\{\sigma_1,\dots,\sigma_q,\tau_1,\dots,\tau_q\}$.

Then, the inclusion map $\iota:\wti X\ra X$ induces isomorphisms on
homology groups with integer coefficients. If, in addition, the spaces
$X$ and $X\sm\{\sigma_1,\dots,\sigma_k,\tau_1,\dots,\tau_k\}$ are
simply connected, for all $1\leq k\leq q$, then there is a~strong
deformation retraction from $X$ to $\wti X$.
\end{thm}
\pr Apply Lemma~\ref{lm:dmt} first to $X$ and pair
$(\sigma_1,\tau_1)$, then to $X\sm\{\sigma_1,\tau_1\}$ and pair
$(\sigma_2,\tau_2)$, etc., until we reach $\wti X$. Take the
concatenation of all the isomorphisms obtained at each step.  
\qed

\subsection{The proof of the main theorem} $\,$

\nin Recall, that by Proposition~\ref{prop:dec}(1), every spherical
word $w$ has a~unique decomposition $w=a_1 v_1 a_1\dots a_t v_t a_t$,
where, for all $1\leq i\leq t$, the word $a_i v a_i$ is circular.

\begin{df}
Let $w=a_1 v_1 a_1\dots a_t v_t a_t$ be a~representation of
a~spherical word, as concatenation of circular ones. We call
$v=a_1^2\dots a_t^2$ the {\bf fundamental subword} of~$w$.
\end{df}

The simplex indexed by the fundamental subword will encode the
topology of $\da_w$.

\begin{df}\label{df:xi}
We define a~function $\xi:\zz_+\ra\zz_+$ as follows.
For an arbitrary nonnegative integer $n$ we set
\[\xi(n):=\begin{cases}
n+1,&\text{ if }n\text{ is even};\\
n-1,&\text{ if }n\text{ is odd}.
\end{cases}\] 
\end{df}
Clearly, the function $\xi$ is a bijection, and $\xi^2$ is the
identity map. The function $\xi$ simply negates the last bit in the
binary representation of a~number, and can also be defined by a~closed
formula $\xi(n)=4\lfloor n/2\rfloor+1-n$.

\begin{df}\label{df:height}
Assume we are given an $n$-tuple $\alpha=(\alpha_1,\dots,\alpha_n)$,
and an index $t$, such that $1\leq t\leq n$, and the numbers
$\alpha_1,\dots,\alpha_{t-1}$ are even.

Let $\beta=(\beta_1,\dots,\beta_n)$ be an $n$-tuple, such that
$\beta\leq\alpha$. We let $h_t(\beta)$ denote the minimal index $k$
between $1$ and $t$, such that $\beta_k\leq\alpha_k-1$. If no such $k$
exists, i.e., if $\beta_i=\alpha_i$, for all $i=1,\dots,t$, then we
set $h_t(\beta):=t$. We call $h_t(\beta)$ the {\bf $t$-height}
of~$\beta$ (w.r.t.\ $\alpha$).

When $v$ is a subword of $w$, we set $h_t(v):=h_t(\beta)$, where
$\beta=(\beta_1,\dots,\beta_t)$, is the left-shifted exponential
presentation of~$v$ (and the $t$-height is taken w.r.t.\ the
exponential presentation of~$w$).
\end{df}

\begin{prop} \label{pr:fund}
Assume a word $w$ is given by its reduced form $a_1^{\alpha_1}\dots a_t^{\alpha_t}$,
such that the numbers $\alpha_1,\dots,\alpha_{t-1}$ are even.
\begin{enumerate}
\item[(1)] If $\alpha_t$ is odd, then the topological space $\da_w$ is contractible.
\item[(2)] If $\alpha_t$ is even, then the $(l(w)-2)$-skeleton of $\da_w$ is
  contractible.  In particular, the topological space $\da_w$ is homotopy equivalent
  to a~$(l(w)-1)$-sphere.
\end{enumerate}
\end{prop}

\pr We will first show the statement~(1). Clearly, in this case, the word $w$ is not
circular. So, by Theorem~\ref{thm:fund} the space $\da_w$ is simply connected. In
order to use the Whitehead's Theorem~\ref{thm:wh1}, we need to see that all integral
homology groups of $\da_w$ vanish.

Given a~subword $v$, for brevity we shall set $h:=h_t(v)$. Let $\Sigma_0$ denote the
set of all left-shifted presentations such that $\beta_h$ is even, and let $\Sigma_1$
denote the set of all left-shifted presentations such that $\beta_h$ is odd. Clearly,
$\Sigma_0$ and $\Sigma_1$ are disjoint, and simplices of $\da_w$ are indexed by
$\Sigma_0\cup\Sigma_1$.

We now define a~matching $\mu_t$ between the sets $\Sigma_0$ and
$\Sigma_1$. Namely, for $v=(\beta_1,\dots,\beta_t)$, we set
\begin{equation}\label{eq:mu}
\mu_t(v):=(\alpha_1,\dots,\alpha_{h-1},\xi(\beta_h), \beta_{h+1},\dots,\beta_t).
\end{equation}

We start by verifying that $\mu_t(v)$ is well-defined. All we need to check is that
$\xi(\beta_h)\leq\alpha_h$. If $\beta_h\leq\alpha_h-1$ this is obvious. Otherwise, we have
$(\beta_1,\dots,\beta_t)=(\alpha_1,\dots,\alpha_t)$. In this case $h=t$, and
$\xi(\alpha_t)=\alpha_t-1$, since $\alpha_t$ is odd.

Next, let us show that $\mu_t(v)$ is a~left-shifted exponential presentation of a~subword
of~$w$. For example, if $w=a^2b^2a$, and $v=(2,1,1)$, then $\mu_t(v)=(2,0,1)$, which is
a~left-shifted exponential presentation of the subword~$a^3$. If $\beta_h$ is even, then
$\xi(\beta_h)=\beta_h+1$, and we obviously get a~left-shifted exponential
presentation. Assume therefore that $\beta_h$ is odd, so $\xi(\beta_h)=\beta_h-1$. Let $m$
be the minimal index, such that $m\geq h+1$, and $\beta_m>0$. If that index does not
exist, we must have $\beta_i=0$, for all $h\leq i\leq t$, and so $\mu_t(v)$ is obviously
left-shifted. Else, the index $m$ is well-defined, and we have
\[\mu_t(v)=(\alpha_1,\dots,\alpha_{h-1},\beta_h-1,\dots,0,\beta_m,\dots).\]  
If $a_m\neq a_h$, then this $t$-tuple is for sure left-shifted. If, on the other hand,
$a_m=a_h$, then the $t$-tuples
\vspace{-1ex}
\begin{alignat*}{2}
\beta  &=(\alpha_1,\dots,\alpha_{h-1},\beta_h,  0,\dots,0,\beta_m,\dots) &\text{\rm and }\\
\beta' &=(\alpha_1,\dots,\alpha_{h-1},\beta_h+1,0,\dots,0,\beta_m-1,\dots)\qquad & 
\end{alignat*} 
are both exponential presentations of $v$, where we recall that here
$\beta_h+1\leq\alpha_h$. Since clearly $\beta'\succl\beta$, we can
conclude that the $t$-tuple $\beta$ was not left-shifted to start
with, yielding a~contradiction with our initial assumptions. We
therefore conclude that $\mu_t(v)$ is a~left-shifted exponential
presentation of a~subword of~$w$.

We can next see that
\begin{equation}\label{eq:hmu}
h=h_t(\mu_t(v)).
\end{equation} 
Clearly, the only way this could fail to be true would be if $\xi(\beta_h)=\alpha_h$.  Of
course, this is impossible if $\alpha_h$ is even. If $\alpha_h$ is odd, then $h=t$, and we
still obtain the identity~\eqref{eq:hmu}.

Furthermore, it follows immediately from $\xi^2=\id$ and \eqref{eq:hmu} that
$\mu_t^2=\id$.  In fact, the map $\mu_t$ provides a~bijection between $\Sigma_0$
and~$\Sigma_1$.

Let us now consider the set of pairs
$\{(\sigma,\mu_t(\sigma))\,|\,\sigma\in\Sigma_0\}$. By what we have
shown, this is a complete decomposition of the set of simplices of
$\da_w$.  Let us now order these pairs in any order which does not
increase the dimension of $\sigma$. A~crucial property which we have
here is the following: if $\sigma\in\Sigma_0$, and $\gamma$ covers
$\sigma$ in $\cp(\da_w)$, then either $\gamma=\mu_t(\sigma)$, or
$\gamma\in\Sigma_0$. To see that the suggested order is actually
a~collapsing order, simply check the three conditions of
Definition~\ref{df:collapsing}. The first two conditions follow from
the construction of $\mu_t$, and the last one follows from the above
mentioned property.  We can thus apply Theorem~\ref{thm:dmt} to
conclude that the space $\da_w$ is contractible.


We shall now show the statement (2), so assume that $\alpha_t$ is even. We can clearly
assume that $\alpha_1+\dots+\alpha_t\geq 4$, since the claim is trivially true for
$w=a_1^2$.  It is now easily seen that the matching $\mu_t$ is still well-defined and that
all the simplices of $\da_w$ are matched, except for the single top-dimensional simplex.

Let $X$ be the complex obtained from $\da_w$ by removing the single top-dimensional
cell. Note, that this cell has dimension~$\geq 3$. The matching $\mu_t$ implies that all
integral homology groups of $X$ vanish. Furthermore, $X$ is simply connected, since it is
obtained from a simply connected space $\da_w$ by removing a cell of dimensional at least
$3$. Whitehead's theorem now implies that $X$ is contractible.  Since the subcomplex
inclusion is a~cofibration, this subcomplex can be shrunk to a~point, yielding a homotopy
equivalence. Attaching the top simplex onto this point yields a~sphere.  \qed

\begin{figure}[hbt]

  \input{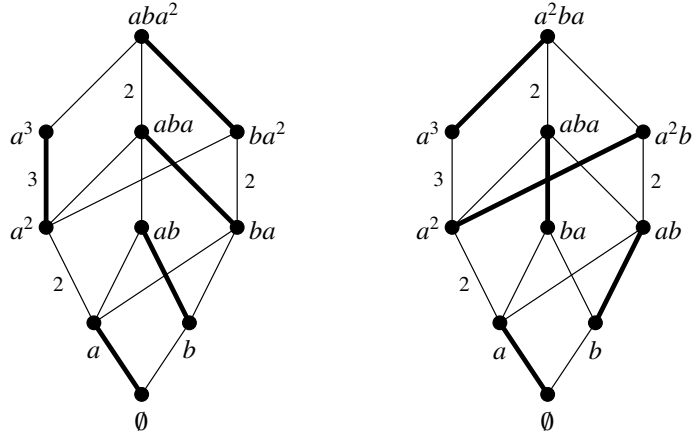}  

\caption{The face posets of $\da_{aba^2}$ and $\da_{a^2ba}$ and 
the induced matchings of simplices.}
\label{fig:aba2}
\end{figure}

\begin{prop}\label{prop:reduc}
Assume a word $w$ is given by its reduced form $a_1^{\alpha_1}\dots
a_t^{\alpha_t}$, such that not all $\alpha_i$ are even. Let $k$ be the
minimal index, such that $\alpha_k$ is odd, and assume that $k\leq
t-1$.  Let $\tilde w$ be obtained from $w$ by removing one letter
$a_{k+1}$ from the power $a_{k+1}^{\alpha_{k+1}}$, i.e., we set
\[\tilde w:=a_1^{\alpha_1}\dots a_k^{\alpha_k} a_{k+1}^{\alpha_{k+1}-1}
a_{k+2}^{\alpha_{k+2}}\dots a_t^{\alpha_t};\] 
note, that when $\alpha_{k+1}=1$ and $a_k=a_{k+2}$ this is not
a~reduced form of the word $\tilde w$.  

Then there exists a strong deformation retraction from $\da_w$ to
$\da_{\tilde w}$.
\end{prop}
\pr The simplices of $\da_w$ which do not belong to $\da_{\tilde w}$ can be
indexed by $(k+1)$-shifted exponential presentations $(\beta_1,\dots,\beta_t)$,
such that $\beta_{k+1}=\alpha_{k+1}$.  By Proposition~\ref{prop:pshift}, such
a~presentation is unique.  The map $\mu_k$ from Proposition~\ref{pr:fund}
provides a complete matching on this set of simplices.  \qed

\vspace{5pt}

\nin
We are now ready to state and to prove our main theorem.

\begin{thm} \label{thm:main}
Let $w$ be an arbitrary word.
\begin{enumerate}
\item[(1)] If $w$ is not spherical, then the $\da$-complex $\da_w$ is
  contractible.
\item[(2)] Assume $w$ is spherical, and let $v=a_1^2\dots a_t^2$ be its fundamental
  subword. Then there exists a strong deformation retraction from $\da_w$ to the
  subcomplex $\da_v$. In particular, as a~topological space $\da_w$ is homotopy equivalent
  to a~$(2t-1)$-sphere.
\end{enumerate} 
\end{thm}

\pr Using Proposition~\ref{prop:reduc} we can reduce every spherical word to the
fundamental subword. If the original word is not spherical. we can reduce it to the word
of the form $a_1^2\dots a_n^2 a$. Flipping the word and applying the same proposition, we
arrive at the word $a$, giving just a point. This means, that the original space was
contractible.  \qed

\vspace{5pt}

Two examples of matchings produced by Theorem~\ref{thm:main} are shown on
Figure~\ref{fig:aba2}. Note, that even though $\da_{aba^2}$ and $\da_{a^2ba}$ are
isomorphic as cell complexes (the indexed word is flipped), the final matchings yielded by
the theorem are quite different.

\subsection{Example of an application: the alternating words} $\,$
\label{ssect:abab}



\nin Recall, that $\alt(n)$ denotes the word $w=ababababa\dots$, such that
$l(w)=n$. Depending on the parity of $n$ we either have $w=(ab)^q$ or $w=(ab)^q a$.

The specific matching given by Theorem~\ref{thm:main} for $w=\alt(n)$
will be 
\begin{align}\label{eq:rul}
(a^2b^2)^kab\sigma&\leftrightarrow(a^2b^2)^kb\sigma, \\
(a^2b^2)^ka^2ba\sigma&\leftrightarrow(a^2b^2)^ka^3\sigma, \notag\\
(a^2b^2)^k&\lra(a^2b^2)^ka,\notag\\
(a^2b^2)^ka^2&\lra(a^2b^2)^ka^2b,\notag
\end{align}
for all $k\geq 0$, and all words $\sigma$.  In particular, $a$ is matched with the empty
simplex. The matching rules will of course only be applied if both simplices are in
$\da_w$. In particular, it is easy to see that when $3$ divides $l(w)$, there will be one
unmatched simplex, namely the one indexed by the fundamental subword of~$w$.

\begin{thm}\label{thm:aba}
When $3$ does not divide $n$, the $\da$-complex $\da_{\alt(n)}$ is collapsible.  If $3$
divides $n$, then one can collapse $\da_{\alt(n)}$ onto $\da_v$, where $v$ is the
fundamental subword of $\alt(n)$. In particular, we have
\begin{equation}
\label{eq:abab}
\da_{\alt(n)}\simeq
\begin{cases}
S^{2n/3-1}, & \textrm{\rm if } 3 \text{ \rm divides } n;\\
\textrm{\rm point,} & \textrm{\rm otherwise.}
\end{cases}
\end{equation}
\end{thm}
\pr The only strengthening of the general theorem here is that the
strong deformation retraction is replaced by collapses in
$\da$-complexes.  This can be done, since all pairs $(\sigma,\tau)$ in
the collapsing order prescribed by Theorem~\ref{thm:main} satisfy
$[\sigma:\tau]=\pm 1$, as can be seen by direct examination of the
rules~\eqref{eq:rul}.  \qed

\vspace{5pt}

\nin So the alternating words for which we get non-trivial topology
are $w=aba$, with $\da_{aba}\simeq S^1$, $w=(ab)^3$, with
$\da_{(ab)^3}\simeq S^3$, $w=(ab)^4a$, with $\da_{(ab)^4a}\simeq S^5$,
etc.

\subsection{Last remarks} $\,$



\nin The following three facts are easy verifications which are left
to the reader.

\begin{itemize}
\item All words for which $\da_w$ is a pseudomanifold have the reduced
  form $a_1^2\dots a_t^2$.

\item All words for which $\da_w$ is a manifold have the reduced form
  $a_1^2\dots a_t^2$, with an additional condition $a_i\neq a_j$, for
  $i\neq j$.

\item When $w=a_1^{\alpha_1}\dots a_t^{\alpha_t}$, such that
  $\alpha_i\geq 2$ for all~$i$, then both the $\da$-complex $\da_w$,
  as well as the simplicial complex $\bd^2\da_w$ are not collapsible.
\end{itemize}

\vspace{5pt} 

\nin In particular, the non-spherical words of the type
$a_1^{\alpha_1}\dots a_t^{\alpha_t}$, with $\alpha_i\geq 2$ for
all~$i$, provide a rich source of contractible, but not collapsible
simplicial complexes. The classical Dunce hat is the special case of
that given by the word $w=a^3$.

\end{document}